# A matheuristic approach for an integrated lot-sizing and scheduling problem with a period-based learning effect

Mohammad Rohaninejad [a], Behdin Vahedi-Nouri [b], Reza Tavakkoli-Moghaddam [b,c], Zdeněk Hanzálek [a,*]

[a] *Czech Institute of Informatics, Robotics, and Cybernetics, Czech Technical University in Prague, Prague, Czech Republic*

[b] *School of Industrial Engineering, College of Engineering, University of Tehran, Tehran, Iran*

[c] *Research Center of Performance and Productivity Analysis, Istinye University, Istanbul, Turkey*

**E-mail addresses:**

Mohammad Rohaninejad: mohammad.rohani.nezhad@cvut.cz

Behdin Vahedi-Nouri: b.vahedi@ut.ac.ir

Reza Tavakkoli-Moghaddam: tavakoli@ut.ac.ir

Zdeněk Hanzálek: zdenek.hanzalek@cvut.cz

[*] Corresponding author.

*E-mail address:* zdenek.hanzalek@cvut.cz (Z. Hanzálek).

**Abstract**

This research investigates a multi-product capacitated lot-sizing and scheduling problem incorporating a novel learning effect, namely the period-based learning effect. This is inspired by a real case in a core analysis laboratory under a job shop setting. Accordingly, a Mixed-Integer Linear Programming (MILP) model is extended based on the big-bucket formulation, optimizing the total tardiness and overtime costs. Given the complexity of the problem, a cutting plane method is employed to simplify the model. Afterward, three matheuristic methods based on the rolling horizon approach are devised, incorporating two lower bounds and a local search heuristic. Furthermore, a post-processing approach is implemented to incorporate lot-streaming possibility. Computational experiments demonstrate: 1) the simplified model performs effectively in terms of both solution quality and computational time; and 2) although the model encounters challenges with large-scale instances, the proposed matheuristic methods achieve satisfactory outcomes; and 3) it can be inferred that the complexity of the models and solution methods are independent of the learning effect; however, the value of learning effect may impact the performance of the lower bounds; 4) in manufacturing settings, where the lot-streaming is possible, incorporating post-processing can drastically improve the objective function; 5) the impact of the period-based learning effect in the results is significant, and the model's sensitivity to time-based parameters (e.g., learning rate) is more than cost-based ones (e.g., tardiness cost).

## 1. Introduction

Under numerous realistic circumstances, the processing time of a repetitive operation involving humans may be reduced due to a phenomenon known as the learning effect [1]. This concept was first studied by Wright [2] in the airplane industry. Later, numerous studies have upheld the existence of the learning effect and its considerable impacts on various industry settings [3, 4]. Moreover, considering the current trend of the market under industry 4.0 and 5.0 revolutions, the manufacturing industry must handle mass-customized and mass-personalized products with short life cycles in a human-centric fashion [5-7]. Therefore, including the learning effect in production scheduling is more important than ever.

Two main learning effects mostly explored in scheduling problems are the *Position-Based Learning Effect* (PBLE) and the *Sum-of-Processing-Time-Based Learning Effect* (SPTBLE). In the PBLE, the learning occurs according to the number of jobs processed before a given job (i.e., setting up machines, deciphering data, and cleaning machines). In contrast, in the SPTBLE, learning happens according to the summation of original processing times of jobs processed before a given time (i.e., welding, painting, and maintenance activities). The incorporation of the learning effect into scheduling problems has broadly been investigated in different environments, including single machine [8, 9], parallel machine [10-12], flow shop [13-15], and job shop [16].

We scrutinize a specific learning effect inspired by a real case in a core analysis laboratory. This laboratory is specialized in examining the geologic characteristics of rock samples. During a year, the laboratory gets involved in several projects with different specified requirements. Based on the internal survey in the laboratory, during consecutive periods of a planning horizon, learning represents itself through workers' error reduction in parameters setting and equipment calibration, increasing measurement accuracy, improving sampling, increasing knowledge and comprehension of reference standards for results validation, enriching mutual understanding with clients, and decreasing conflicts and reworks. Due to these facts, the impact of the learning effect on processing time reduction is quite evident, while this reduction in some operations varies between 15 and 40 percent.

Although the workers' learning has experimentally been accepted in the laboratory processes, it cannot be modeled by classical PBLE and SPTBLE, which are mostly suitable for short-term scheduling. The main reason is that in the short term (during a single period), the processing time of a repetitive operation is not necessarily decreased. It does not follow a specified pattern since the characteristics of samples are different according to their porosity, permeability, density, the type and volume of their fluids, and so on. Therefore, a new kind of learning effect, named the *period-based learning effect*, is formulated for the medium-term horizon based on the laboratory's experimental observations. In this formulation, the actual time needed to process a certain job in a given period decreases according to the number of elapsed periods within the planning horizon. This type of learning is a modification of the PBLE formulation. It is assumed that if job $j$ is processed during period $t$, its actual processing time is equal to $\left(p_j t^a\right)$, in which $p_j$ is job $j$'s original processing time without any learning effect and $a < 0$ is the learning index. The period-based learning effect is also applicable in several cases in manufacturing settings where the workers' learning occurs in medium- or long-term horizon, including in industries requiring high-skill craftsmanship (e.g., assembly in the aerospace industry), implementing continuous improvement programs (e.g., lean manufacturing), practicing knowledge transferring over regular periods, to name a few.

According to the characteristics of our problem, it is defined as an Integrated Lot-sizing and Scheduling Problem (ILSP). A few researchers have attempted to incorporate the learning phenomenon in pure lot-sizing problems [17-19], which are generally modeled based on the size of a lot so that a larger lot leads to a shorter processing time. As far as the authors know, an attempt yet to be made in the literature to consider the learning effect in ILSPs. Table 1 summarizes these studies and highlights the contributions of this study.

**Table 1.** Summary of studies considering learning effect in lot-sizing

| References | Lot-sizing | Scheduling | Multi-product | Environment | Learning type | Objective function | Solution method |
|---|---|---|---|---|---|---|---|
| Jaber, Bonney and Moualek [17] | ✓ | | | Single machine | PBLE | Sum of production and holding cost | Mathematical modeling |
| Chen and Tsao [18] | ✓ | | ✓ | Single machine | PBLE | Profit | iterative search algorithm |
| Jeang and Rahim [19] | ✓ | | | Single machine | PBLE | Sum of inventory, production and setup costs | Mathematical modeling |
| Jaber and Peltokorpi [20] | | | | Single machine | Incorporating Worker number in PBLE | Sum of inventory, worker and setup costs | Mathematical modeling |
| Our study | ✓ | ✓ | ✓ | Job shop | Period-based learning effect | Sum of overtime and tardiness costs | MILP modeling & Matheuristics |

ILSPs are among Non-deterministic Polynomial hard (NP-hard) problems [20]. So far, several solution methods have been presented to overcome its complex nature. A number of them are founded on mathematical models, including the rolling horizon method [21-24], fix-and-relax heuristics [25-27], fix-and-optimize heuristics [28-31], decomposition heuristic [32], and Lagrangian-based heuristic approach [33]. Other methods, including simulated annealing [34], branch-and-check [35], artificial bee colony algorithm [36], hybrid meta-heuristic algorithm [37, 38], fruit fly algorithm [39], hyper-heuristic [40], genetic algorithm [41], and satisfiability modulo theories (SMT) formulation [42] have also been developed for this purpose.

As demonstrated in the literature, the rolling horizon heuristic is an efficient method that has mainly been employed for dynamic ILSPs, where the amount of demand, resources, system costs, and other parameters of the production plan are uncertain and revealed during the planning horizon. Nevertheless, this heuristic may also apply to a static environment where the quantity and timing of products' demand are specified at the beginning of the planning horizon. This paper develops three versions of the static rolling horizon procedure based on a new big-bucket time model by applying two lower bounds and a local search. The aim is to find an optimal solution for an ILSP in a job shop setting (ILSP-JSP) that minimizes the total cost of tardiness and overtime.

To sum up, the core contributions of this study can be summarized below.

I. We incorporate a new learning effect into ILSP-JSP inspired by a real case in a core analysis laboratory.

II. We develop a big-bucket time mathematical model for the given problem.

III. We simplify the proposed mathematical model using a cutting plane method.

IV. Since the problem is inherently complicated, we devise three matheuristic methods based on the rolling horizon approach incorporating two lower bounds and a local search method to achieve high-quality solutions within a reasonable time.

V. We carry out computational experiments to weigh the performance of the developed models and matheuristic methods.

The rest of this study is arranged as follows. The ILSP-JSP and its assumptions are elaborated in Section 2. Section 3 presents the corresponding big-bucket Mixed-Integer Linear Programming (MILP) model. Section 4 illustrates the simplified model founded on a cutting plane. Section 5 explains the MILP-based heuristic. Section 6 introduces a post-processing approach to integrate lot-streaming into the final model outputs, wherever applicable. Section 7 provides the results of the computational experiments. Ultimately, Section 8 provides the conclusion and a few directions for future studies.

## 2. Problem description

As stated earlier, the main focus of this study is to explore ILSP-JSP with the learning effect. The planning horizon contains a few periods, each with definite time units. In addition, the capacity of each machine during each period is given and may differ from one period to another. Moreover, this capacity can be expanded by overtime, which incurs an overtime cost. The capacity of a machine under regular and overtime conditions in every period is defined in advance. Additionally, the following assumptions are considered.

- The structure of the bill of material is supposed to be “series type”, and the consumption coefficient from one level to another is set to 1.
- Time horizon is finite with $T$ periods.
- Machines are available all the time.
- Machines can, at the most, process one job at a time.
- Each operation can be performed at most once during a period.
- Because of technical limitations, the working time of a machine in planning periods is limited.
- Despite the limited capacity of machines, it is possible to carry out overtime on machines by paying the overtime costs.
- Back-orders are permitted.
- Demands, processing times, overtime costs, learning index, and due dates are deterministic within the planning periods.
- Customers’ demands only exist for the final operation of a job.
- In-process inventory and shortage are permitted; jobs may wait at a machine to process its next operation until the machine becomes idle.
- Job processing may be continued without any interruption.
- Demand satisfaction is indivisible, meaning all units processed to fulfill a specific demand must be delivered together as a single batch.
- Learning is presumed to occur between two successive periods.
- No lead time is assumed.
- The total tardiness and overtime cost are considered as the objective function.

## 3. Big-bucket time formulation

This section presents an MILP model based on a big-bucket time pattern with $J$ jobs and $M$ machines. This model is derived from the big-bucket time model presented by Rohaninejad et al. [37]. In this model, each job includes an order of operations $O_{j,h}$ ($h = 1, \ldots, h_j$), where $O_{j,h}$ denotes the $h$-th operation regarding job *j.* The associated indices, parameters, variables, and proposed MILP model are as follows.

**Indices:**

| | |
|---|---|
| $j, k$ | Index of jobs ($j, k = 1, .., J$) |
| $h, l$ | Index of operations ($h, l = 1, .., h_j$) |
| $r, r'$ | Index of positions ($r, r' = 1, .., N$) |
| $m, m'$ | Index of machines ($m, m' = 1, .., M$) |
| $t, t', t''$ | Index of periods ($t, t', t'' = 1, .., T$) |

**Parameters:**

| | |
|---|---|
| $TC_j$ | Tardiness cost related to job *j* |
| $OC_m$ | Overtime cost on machine *m* per time unit |
| $C_{m,t}$ | Machine *m's* capacity during period *t* |
| $p_{j,h}$ | Processing time required for producing one item of operation $O_{j,h}$ |
| $a_{j,h}$ | Learning index for operation $O_{j,h}$ |
| $d_{j,t}$ | Due date for job *j* during period *t* |
| $D_{j,t}$ | Demand for job *j* during period *t* |
| $O_{m,t}$ | Maximum allowable overtime regarding machine *m* during period *t* |
| $R$ | Allowable number of positions determined for a machine during a period |
| $h_j$ | Number of operations needed to complete job *j* and the $h$-th operation of job *j* |
| $m(j,h)$ | Machine, which is supposed to process operation $O_{j,h}$ |
| $L$ | Length of period *t* |
| $G$ | A large number |

**Variables:**

| | |
|---|---|
| $I_{j,h}^{t}$ | Amount of inventory regarding operation $O_{j,h}$ at the end of period *t* |
| $\varphi_j^t$ | Shortage regarding job *j* at the end of period *t* ($\varphi_j^T = 0$) |
| $X_{j,h}^{t,r}$ | Quantity of produced items regarding operation $O_{j,h}$ in position *r* on its eligible machine during period *t* (the lot size). |
| $o_{m,t}$ | Amount of overtime processing for machine *m* during period *t* |
| $s^{t,m.r}$ | Start time regarding the *r*-th position, on machine *m* during period *t* |
| $f^{t,m.r}$ | Finish time regarding the *r*-th position on machine *m* during period *t* |

$\eta_{j,h}^{t,r,r'}$ — Upper bound of $X_{j,h}^{t,r'}$ if the $r'$-th position on machine $m'$ starts after the start of the $r$-th position on machine $m$ in period $t$ where $m' = m(j,h)$ and $m = m(j,h+1)$

$\omega_{j,h}^{t,r,r'}$ — Upper bound of $X_{j,h}^{t,r'}$ if the $r'$-th position on machine $m'$ starts before the start of the $r$-th position on machine $m$ in period $t$ where $m' = m(j,h)$ and $m = m(j,h+1)$

$F_{j,t}$ — Finish time of job $j$ during period $t$

$i_{j,h}^{t,r}$ — Inventory lower bound regarding $O_{j,h}$ before the start of position $r$, on machine $m$ during period $t$ where $m = m(j,h+1)$

$Tr_{j,t}$ — Tardiness of job $j$ at period $t$

$q^{t,m.r,m',r'}$ — 1 if the $r'$-th position on machine $m'$ starts after the start of the $r$-th position on machine $m$ in period t; 0, otherwise

$y_{j,h}^{t,r}$ — 1 if operation $O_{j,h}$ is set up in position $r$ on its eligible machine during period $t$; 0, otherwise

$\rho_j^{t,t'}$ — 1 if $D_{j,t}$ is satisfied in period $t'$ $(t' \geq t)$; 0, otherwise

Now, we can introduce the following model.

**Model-I:**

$$\text{Min} \sum_{j=1}^{J}\sum_{t=1}^{T} TC_j.Tr_{j,t} + \sum_{m=1}^{M}\sum_{t=1}^{T} o_{m,t}.OC_m \tag{1}$$

s.t.

$$I_{j,h_j}^{t-1} + \sum_{r=1}^{R} X_{j,h_j}^{t,r} - I_{j,h_j}^{t} + \varphi_j^{t} - \varphi_j^{t-1} = D_{j,t} \qquad j = 1,\dots,J\,;\, t = 1,\dots,T; \tag{2}$$

$$I_{j,h}^{t-1} + \sum_{r=1}^{R} X_{j,h}^{t,r} - I_{j,h}^{t} = \sum_{r=1}^{R} X_{j,h+1}^{t,r} \qquad j = 1,\dots,J\,;\; t = 1,\dots,T; h = 1,\dots h_j - 1 \tag{3}$$

$$\sum_{j=1}^{J}\sum_{(h=1|m(j,h)=m)}^{h_j}\sum_{r=1}^{R} X_{j,h}^{t,r}.p_{j,h}.t^{a_{j,h}} \leq C_{m,t} + o_{m,t} \qquad t = 1,\dots,T;\; m = 1,\dots,M \tag{4}$$

$$X_{j,h}^{t,r} - G.y_{j,h}^{t,r} \leq 0 \qquad j = 1,\dots,J\,;\; t = 1,\dots,T;\; h = 1,\dots h_j;\; r = 1,\dots,R \tag{5}$$

$$\sum_{j=1}^{J}\sum_{(h=1|m(j,h)=m)}^{h_j} y_{j,h}^{t,r} \leq 1 \qquad t = 1,\dots,T; r = 1,\dots,R\,;\; m = 1,\dots,M \tag{6}$$

$$\sum_{t''=1}^{t}\sum_{r=1}^{R} X_{j,h_j}^{t'',r} \geq \sum_{t''=1}^{t}\sum_{t'=t''}^{t} \rho_j^{t'',t'}.D_{j,t''} \qquad j = 1,\dots,J\,; t = 1,\dots,T \tag{7}$$

$$\sum_{r'=1}^{R}\left(\eta_{j,h}^{t,r,r'} + \omega_{j,h}^{t,r,r'}\right) - \sum_{r'=r}^{R} X_{j,h+1}^{t,r'} = I_{j,h}^{t} - i_{j,h}^{t,r} \qquad j = 1,\dots,J\,;\; h = 1,\dots h_j - 1\,;\; t = 1,\dots,T\,;\; r = 1,\dots,R \tag{8}$$

$$i_{j,h}^{t,r} \geq X_{j,h+1}^{t,r} \qquad j = 1, \dots, J\,;\ t = 1, \dots, T\,;\ h = 1, \dots h_j - 1; r = 1, \dots, R \tag{9}$$

$$X_{j,h}^{t,r'} \leq \eta_{j,h}^{t,r,r'} + \left(3 - q^{t,m,r,m',r'} - y_{j,h}^{t,r'} - y_{j,h+1}^{t,r}\right).G \qquad j = 1, \dots, J\,;\ h = 1, \dots h_j - 1\,;\ t = 1, \dots, T\,;\ r\&r' = 1, \dots, R\,;\ m' = m(j,h); m = m(j,h+1) \tag{10}$$

$$X_{j,h}^{t,r'} - \frac{\left(s^{t,m.r} - s^{t,m'r'}\right)}{p_{j,h}.t^{a_{j,h}}} \leq \omega_{j,h}^{t,r,r'} + \left(2 + q^{t,m,r,m',r'} - y_{j,h}^{t,r'} - y_{j,h+1}^{t,r}\right).G \qquad j = 1, \dots, J\,;\ h = 1, \dots h_j - 1\,;\ t = 1, \dots, T\,;\ r\&r' = 1, \dots, R\,;\ m' = m(j,h); m = m(j,h+1) \tag{11}$$

$$\left(s^{t,m.r} - s^{t,m'r'}\right) - \left(1 - \mathrm{q}^{t,m,r,m'r'}\right).\mathrm{G} \leq 0 \tag{12}$$

$$\left(s^{t,m'r'} - s^{t,m.r}\right) - \mathrm{q}^{t,m,r,m'r'}.\mathrm{G} \leq 0 \tag{13}$$

$t = 1, \dots, T\,;\ r\&r' = 1, \dots, R\,;\ m\&m' = 1, \dots, M$

$$\eta_{j,h}^{t,r,r'} - q^{t,m.r,m'r'}G \leq 0 \qquad j = 1, \dots, J\,;\ h = 1, \dots h_j - 1\,;\ t = 1, \dots, T\,;\ r\&r' = 1, \dots, R\,;\ m' = m(j,h); m = m(j,h+1) \tag{14}$$

$$\omega_{j,h}^{t,r,r'} - \left(1 - q^{t,m,r,m'r'}\right).G \leq 0 \qquad j = 1, \dots, J\,;\ h = 1, \dots h_j - 1\,;\ t = 1, \dots, T\,;\ r\&r' = 1, \dots, R\,;\ m' = m(j,h); m = m(j,h+1) \tag{15}$$

$$f^{t,m.r} = s^{t,m.r} + \sum_{j=1}^{J} \sum_{(h=1|m(j,h)=m)}^{h_j} X_{j,h}^{t,r}.p_{j,h}.t^{a_{j,h}} \qquad t = 1, \dots, T\,; m = 1, \dots, M;\ r = 1, \dots, R \tag{16}$$

$$s^{t,m.r+1} \geq f^{t,m.r} \qquad t = 1, \dots, T\,; m = 1, \dots, M;\ r = 1, \dots, R \tag{17}$$

$$f^{t,m.r} \leq L.t \qquad t = 1, \dots, T\,; m = 1, \dots, M;\ r = 1, \dots, R \tag{18}$$

$$s^{t,m.r} \geq L.(t-1) \qquad t = 1, \dots, T\,; m = 1, \dots, M;\ r = 1, \dots, R \tag{19}$$

$$F_{j,t} \geq f^{t',m.r} - \left(2 - \rho_j^{t,t'} - y_{j,h_j}^{t',r}\right).G \qquad j = 1, \dots, J\,;\ t = 1, \dots, T\,;\ t' = t, \dots, T\,; r = 1, \dots, R \tag{20}$$

$$Tr_{j,t} \geq F_{j,t} - d_{j,t} \qquad j = 1, \dots, J\,;\ t = 1, \dots, T \tag{21}$$

$$\sum_{r=1}^{R} y_{j,h}^{t,r} \leq 1 \qquad j = 1, \dots, J\,;\ h = 1, \dots h_j\,;\ t = 1, \dots, T \tag{22}$$

$$o_{m,t} \leq O_{m,t} \qquad t = 1, \dots, T\,; m = 1, \dots, M \tag{23}$$

$$\sum_{t'=t}^{T} \rho_j^{t,t'} = 1 \qquad j = 1, \dots, J\,;\ t = 1, \dots, T \tag{24}$$

Equation (1) expresses the objective function comprising the sum of the tardiness and overtime costs. Constraint (2) ensures inventory balance for the final operation of each job across macro periods, while Constraint (3) maintains this balance for the remaining operations. Constraint (4) observe machine capacity, preventing the overutilization of resources. Constraints (5) and (6) stipulate the position of operation $O_{j,h}$ on machine $m$, ensuring that only one operation can be performed in a given position. Constraint (7) ensures that the cumulative production of job $j$ by the end of period $t$ is at least equal to the cumulative demand for job $j$ that has been fulfilled by the end of period $t$. Constraint (8) ensures the inventory balance inside a macro period, and Constraint (9) defines the minimum quantity of inventory needed to produce an item in the position $r$. These constraints ensure that sufficient inventory of $O_{j,h}$ is available before starting $O_{j,h+1}$ within a given period $t$.

Constraints (10) and (11) determine the upper bounds for producing an item in a position by respecting the inventory balance at the associate period. These constraints prevent overestimating the variable $i_{j,h}^{t,m,r}$, which could otherwise result in infeasible solutions. Meanwhile, the combination of these constraints and Constraints (8) and (9) allows $O_{j,h+1}$ to overlap with $O_{j,h}$ or even start before $O_{j,h}$, provided there is sufficient inventory of $O_{j,h}$ available at the beginning of $O_{j,h+1}$. While determining the inventory balance in a macro period, the production scheduling in that period must be known. Therefore, Constraints (12) - (15) are applied to specify the transposition or simultaneity of two different production positions within a period. Constraints (16) - (19) specify the start time and finish time of each position. Constraints (20) and (21) calculate the tardiness regarding job $j$. Constraint (22) ensures that an operation can be done once within a period. Constraint (23) specifies the ranges of overtime variables. Constraint (24) guarantees that the demand $D_{j,t}$ must be fulfilled in a period $t'$, where $t' \geq t$.

## 4. Cutting plane

Model-I in Section 3 is not capable of tackling large-scale problems. Therefore, this model shall be simplified. This section obtains an upper bound for the model using a cutting plane described below.

***Lemma 1.*** The optimal objective value for Model-I is bounded on top by adding Relation (25).

$$s^{t,m.r} \geq f^{t,m',r'} + \left(y_{j,h}^{t,r'} + y_{j,h+1}^{t,r} - 2\right).G \qquad \begin{array}{l} j = 1,\ldots,J\,;\ h = 1,\ldots h_j - 1\,; \\ t = 1,\ldots,T\,;\ r\&r' = 1,\ldots,R \end{array} \qquad (25)$$

**Proof:** To prove this lemma, it is important to note that Relation (25) ensures that if two consecutive operations of a job are processed within the same period $t$, the predecessor operation must be completed before the successor operation begins. Assume that $O_{j,h}$ is executed in the $r'$-th position on machine $m'$, and $O_{j,h+1}$ is executed in the $r$-th position of machine $m$. By applying Relation (25), the inventory of $O_{j,h}$ available at the start of the $r$-th position on machine $m$ as $i_{j,h}^{t,r}$ is given by:

$$i_{j,h}^{t,r} = I_{j,h}^{t-1} + X_{j,h}^{t,r'} \tag{26}$$

Furthermore, by incorporating Constraint (3) into Equation (26), the inventory $i_{j,h}^{t,r}$ can be expressed as:

$$i_{j,h}^{t,r} = I_{j,h}^{t} + X_{j,h+1}^{t,r} - X_{j,h}^{t,r'} + X_{j,h}^{t,r'} = I_{j,h}^{t} + X_{j,h+1}^{t,r} \tag{27}$$

On the other hand, based on Constraint (8) we have:

$$i_{j,h}^{t,r} = I_{j,h}^{t} + X_{j,h+1}^{t,r} - \left(\eta_{j,h}^{t,r,r'} + \omega_{j,h}^{t,r,r'}\right) \tag{28}$$

By comparing Equations (27) and (28), we conclude that $\eta_{j,h}^{t,r,r'} + \omega_{j,h}^{t,r,r'} = 0$ when $s^{t,m.r} \geq f^{t,m',r'}$. According to Constraint (9), $i_{j,h}^{t,r} \geq X_{j,h+1}^{t,r}$, which implies that it is permissible to reduce $i_{j,h}^{t,r}$ if it is strictly greater than $X_{j,h+1}^{t,r}$ (i.e., $i_{j,h}^{t,r} > X_{j,h+1}^{t,r}$).

To achieve this reduction in $i_{j,h}^{t,r}$, the term $\eta_{j,h}^{t,r,r'} + \omega_{j,h}^{t,r,r'}$ must take a positive value (i.e., $\eta_{j,h}^{t,r,r'} + \omega_{j,h}^{t,r,r'} > 0$). However, based on Constraints (10) to (15), a positive value of $\eta_{j,h}^{t,r,r'} + \omega_{j,h}^{t,r,r'}$ requires shifting the $r$-th position on machine $m$ to the right relative to the $r'$-th position on machine $m'$. Fig. 1 schematically illustrates how shifting the $r$-th position on machine $m$ to the right results in a positive value for $\eta_{j,h}^{t,r,r'} + \omega_{j,h}^{t,r,r'}$. This figure also shows how varying values of $\eta_{j,h}^{t,r,r'}$ and $\omega_{j,h}^{t,r,r'}$ ensure compliance with the material balance constraint within a single period.

Importantly, such a shift does not increase tardiness or overtime cost. Consequently, Relation (25) consistently provides an upper bound for Model-I, completing the proof. ☐

The performance of this cutting plane in solving the problem is compelling. The computational results show that the upper bounds obtained from this inequality are so close to an optimal solution of Model-I. Its main advantage is the ability to simplify the original model by new constraints without any need to provide a pattern, as is the case with the decomposition or relaxation of the model. Therefore, Model-I is simplified based on the above cutting plane. The new variables and constraints of the simplified MILP model (Model-II) are as follows.

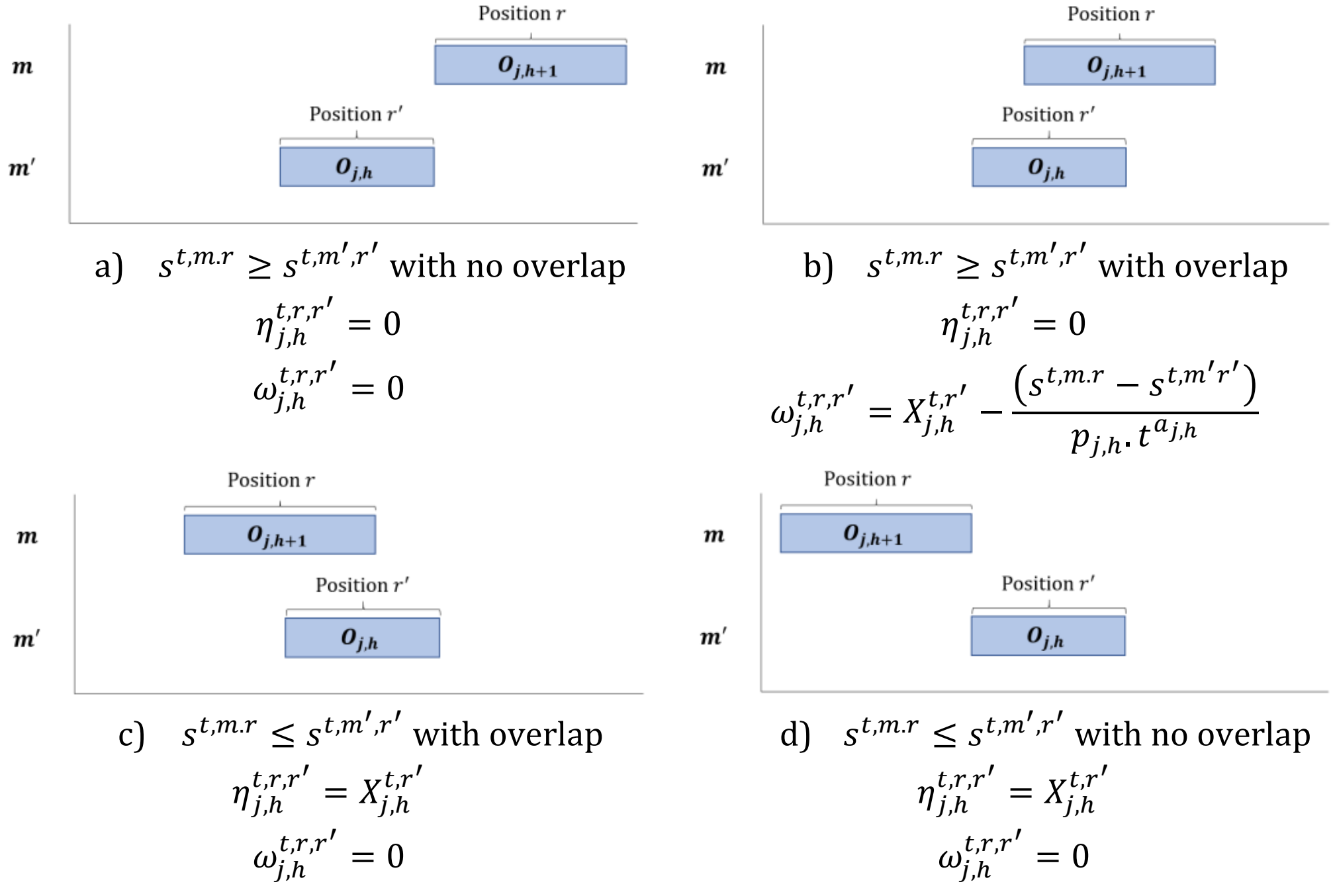

**Fig. 1.** Illustration of the material balance constraint for two consecutive operations within a macro-period

**New variables of Model-II:**

- $z_{j,h,k,l}^{t}$ — 1 if operation $O_{j,h}$ is performed within period $t$ before operation $O_{k,l}$ and $m(j,h) = m(k,l)$; 0, otherwise.
- $x_{j,h}^{t}$ — Processing time regarding operation $O_{j,h}$ within period *t.*
- $s_{j,h}^{t}$ — Operation $O_{j,h}$'s start time during period *t*.
- $f_{j,h}^{t}$ — Operation $O_{j,h}$'s finish time during period *t*.

**Model-II:**

Constraints (21), (23), and (24)

$$x_{j,h+1}^{t} \leq \sum_{t'=1}^{t} x_{j,h}^{t'} - \sum_{t'=1}^{t-1} x_{j,h+1}^{t'} \qquad j = 1,\ldots,J;\ t = 1,\ldots,T;\ h = 1,\ldots,h_j - 1 \tag{29}$$

$$\sum_{t=1}^{T} x_{j,h}^{t} = \sum_{t=1}^{T} D_{j,t} \qquad j = 1,\ldots,J;\ h = 1,\ldots,h_j \tag{30}$$

$$\sum_{t''=1}^{t} x_{j,h_j}^{t''} \geq \sum_{t''=1}^{t} \sum_{t'=t''}^{t} \rho_j^{t'',t'} D_{j,t''} \qquad j = 1,\ldots,J;\ t = 1,\ldots,T \tag{31}$$

$$\sum_{j=1}^{J} \sum_{(h=1|m(j,h)=m)}^{h_j} x_{j,h}^{t}.p_{j,h}.t^{a_{j,h}} \leq C_{m,t} + o_{m,t} \qquad m = 1,\ldots,M;\ t = 1,\ldots,T \tag{32}$$

$$z_{j,h,k,l}^{t} + z_{k,l,j,h}^{t} = 1 \qquad j,k = 1,\ldots,J;\ t = 1,\ldots,T;\ h,l = 1,\ldots,h_j \ |O_{j,h} \neq O_{k,l}\ \&\ m(j,h) = m(k,l) \tag{33}$$

$$f_{j,h}^{t} \leq s_{k,l}^{t} + \left(1 - z_{j,h,k,l}^{t}\right).G \qquad \begin{array}{l} j,k = 1,..,J;\ t = 1,..,T \\ h,l = 1,..,h_j \\ |O_{j,h} \neq O_{k,l}\ \&\ m(j,h) = m(k,l) \end{array} \tag{34}$$

$$s_{j,h+1}^{t} \geq f_{j,h}^{t} \qquad \begin{array}{l} j = 1,..,J;\ t = 1,..,T \\ h = 1,..,h_j - 1 \end{array} \tag{35}$$

$$f_{j,h}^{t} = s_{j,h}^{t} + x_{j,h}^{t}.p_{j,h}.t^{a_{j,h}} \qquad \begin{array}{l} j = 1,..,J;\ t = 1,..,T \\ h = 1,..,h_j \end{array} \tag{36}$$

$$s_{j,h}^{t} \geq L.(t-1) \qquad \begin{array}{l} j = 1,..,J;\ t = 1,..,T \\ h = 1,..,h_j \end{array} \tag{37}$$

$$f_{j,h}^{t} \leq L.t \qquad \begin{array}{l} j = 1,..,J;\ t = 1,..,T \\ h = 1,..,h_j \end{array} \tag{38}$$

$$F_{j,t} \geq f_{j,h_j}^{t'} - \left(1 - \rho_j^{t,t'}\right).G \qquad \begin{array}{l} j = 1,..,J;\ t = 1,..,T \\ t' = t,...,T \end{array} \tag{39}$$

Constraints (29) and (30) determine the relationship between inventory levels, planned quantities, and demands. In other words, these constraints are replaced with Constraints (2) and (3) in the Model-I. Constraint (31) is equivalent to the Constraint (7). Constraint (32) and is substituted for Constraints (4). Constraints (33) and (34) determine the sequence of operations $O_{j,h}$ and $O_{k,l}$ within period $t$. Constraint (35) ensures that the start time of operation $O_{j,h+1}$ is after finishing operation $O_{j,h}$ within period $t$. By adding this constraint, Constraints (8) to (15) are relaxed in the Model-I. This constraint and Constraint (29) guarantee the inventory balances within a macro period. Constraints (36) to (38) are substituted for constraints (16) to (19). Constraint (20) is replaced with Constraint (39). Constraints (21), (23), and (24) are identically repeated in the new model.

First, the complexity of the models to formulate same-sized problems is analyzed to compare Model-I with Model-II. For this purpose, four instances of different sizes (based on $J$, $M$, and $T$) are generated. Table 2 shows the number of constraints, continuous variables, and binary variables, which each model requires to formulate the different instances. For each instance, the total number of operations equals 3×$J$ (every job consists of three operations), and factor $R$ in Model-I equals 3×$J/M$. As it is clear from the results, Model-II needs fewer constraints and variables than Model-I. Hence, Model-II is expected to have a lower complexity than Model-I.

**Table 2.** Comparing the complexity of the two models

| Instance size | | | No. of constraints | | No. of continuous variables | | No. of binary variables | |
|---|---|---|---|---|---|---|---|---|
| $J$ | $M$ | $T$ | Model-I | Model-II | Model-I | Model-II | Model-I | Model-II |
| 5 | 3 | 3 | 4953 | 663 | 2829 | 174 | 705 | 210 |
| 20 | 6 | 5 | 112460 | 7420 | 66830 | 1130 | 18300 | 3000 |
| 50 | 10 | 8 | 1000360 | 42110 | 579280 | 4480 | 181800 | 18600 |
| 100 | 15 | 10 | 4194300 | 136100 | 2528150 | 11150 | 905500 | 62500 |

## 5. Matheuristic methods

Real-world problems are very complex. Thus, employing heuristic and metaheuristic methods is a prevalent way to tackle them [43]. MILP-based heuristics (e.g., rolling horizon method) effectively discover decent solutions for large-sized problems [44], particularly for ILSPs. These heuristic methods are regularly founded on decomposing the original problem into subsets that can be solved simply by an MILP in an iterative way.

This paper develops three distinct rolling horizon heuristics based on two new lower bounds and a local search heuristic, which are established by applying Model-II. A rolling horizon heuristic method is one of the most effective matheuristic methods, and it has mainly been used in dynamic environments with lot-sizing and scheduling problems. Nonetheless, this approach is also suitable for static models, in which the quantity and timing of products' demand are specified during the planning horizon. The original problem is segmented into a series of sub-problems, focusing on earlier and later periods. Accordingly, detailed planning is applied to the earlier periods, while broader, high-level planning is used for the following periods. The size and number of these sub-problems significantly influence the excellence of the ultimate solution and the problem complexity.

The proposed heuristic method segments the planning horizon into $T$ distinct time intervals (i.e., sub-problems), with each interval representing a specific period. Fig. 2 relates how to carry out the proposed heuristic method in an iterative way.

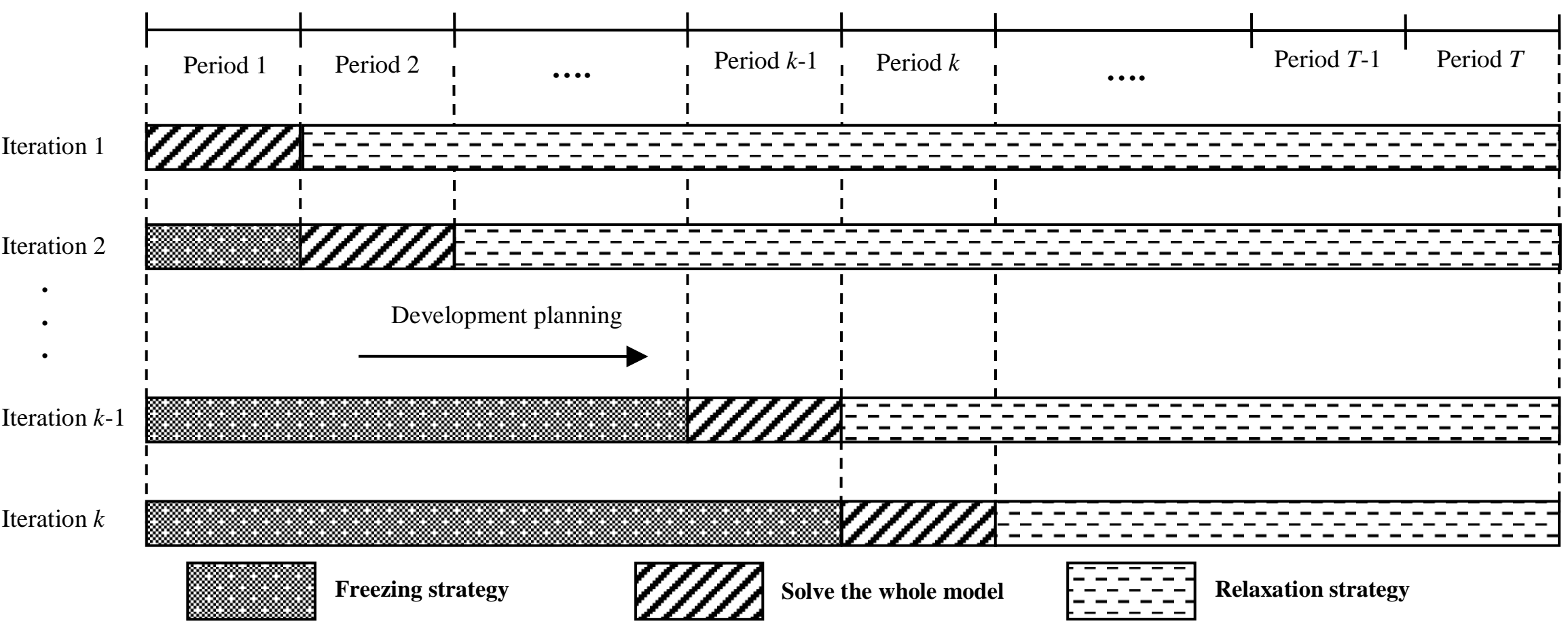

**Fig. 2.** Schematic of the developed rolling horizon method.

At each iteration, three tactics containing the "freezing", "solving the whole model" and "relaxation" can be executed. In the first iteration, the model associated with the first time interval (i.e., period 1) is solved, while for subsequent intervals (i.e., period 2 to $T$), the relaxation tactic simplifies the model by relaxing the binary variables ($y^{t}_{j,h}$ and $z^{t}_{j,h,k,l}$). Likewise, the freezing tactic is employed for sub-problem $k$-1 within the $k$-th iteration. Regarding the freezing tactic, the binary variables associated with period $k$-1 ($y^{k-1}_{j,h}$ and $z^{k-1}_{j,h,k,l}$) are fixed according to the solution obtained by solving the entire model in iteration $k$-1. In addition, the model relevant to the time interval $k$ is regarded as a whole, and the relaxation tactic simplifies the models of sub-problem $k$+1 to $T$.

### *5.1. Lower Bounds*

One drawback of the relaxation strategy is due to the abandonment of control constraints, such as priority constraints and constraints related to the sequence of operations. Clearly, a great distance concerning the conditions between the relaxed problems and the original problem leads to a deterioration of the solution quality. Hence, it is essential to ensure that the conditions of relaxed issues are close to those of the original problem, with a minimum increase in the complexity of the relaxed one. Here, this issue is unraveled by providing two lower bounds.

***Lemma 2. (Lower bound 1).*** The optimal value of Model-II is bounded by Relation (40). In other words, a lower bound is created for Model-II by replacing Constraint (21) in Model-II with Constraint (40).

$$Tr_{j,t} \geq \left( \left( (L(t-1)+1)u_j^t - 1 \right) G + \sum_{h=1}^{h_j} x_{j,h}^t . p_{j,h} . t^{a_{j,h}} \right) - d_{j,t} \qquad t = 1,..,T\,; \; j = 1,..,J \qquad (40)$$

In this constraint, $u_j^t$ is a binary variable equal to 1 if at least one operation of job $j$ is performed within period $t$; otherwise, it takes 0. This variable can be obtained by:

$$\left( \sum_{h=1}^{h_j} x_{j,h}^t \right) - u_j^t . G \leq 0 \qquad j = 1,..,J;\, t = 1,..,T \qquad (41)$$

***Proof.*** The following two cases exist when analyzing Relation (40).

*Case 1.* If there is no operation $O_{k,l}$ $(k \neq j)$ that $z_{k,l,j,h}^t = 1$ for any operations $O_{j,h}$ $(h \leq h_j)$, then $F_{j,t} = L(t-1) + \sum_{h=1}^{h_j} x_{j,h}^t . p_{j,h} . t^{a_{j,h}}$ . In other words, $Tr_{j,t}$ obtained by Relation (40) equals $Tr_{j,t}$ obtained by Relation (21).

*Case 2.* If there is at least one operation $O_{k,l}$ $(k \neq j)$ that $z_{k,l,j,h}^t = 1$ for operation $O_{j,h}$ $(h \leq h_j)$, then $F_{j,t} > L(t-1) + \sum_{h=1}^{h_j} x_{j,h}^t . p_{j,h} . t^{a_{j,h}}$. In other words, $Tr_{j,t}$ obtained by Relation (40) is smaller than $Tr_{j,t}$ obtained by Relation (21).

Combining these two cases shows that $Tr_{j,t}$ obtained by Relation (40) is always smaller than or equal to $Tr_{j,t}$ obtained by Relation (21), and the proof is complete. ☐

***Lemma 3. (Lower bound 2).*** The optimal value of the Model-II is bounded by:

$$Tr_{j,t} \geq \left( \sum_{h=1}^{h_j} ( \sum_{t'=1}^{t} D_{j,t'} - \sum_{t'=1}^{t-1} x_{j,h}^{t'} ) . p_{j,h} . T^{a_{j,h}} \right) - \left( d_{j,t} - \left( L.(t-1) \right) \right) \qquad j = 1,..,J; \; t \in tp_j \qquad (42)$$

In the relaxation strategy, in the relaxed sub-problems, Constraint (21) is substituted with Constraint (42), where $tp_j$ is the period in which the due dates of job *j* ($d_{j,t}$) take place.

***Proof.*** The remaining time for job $j$ within period $t$ equals $\left( d_{j,t} - \left( L.(t-1) \right) \right)$ and the remaining work equals $\left( \sum_{h=1}^{h_j} ( \sum_{t'=1}^{t} D_{j,t'} - \sum_{t'=1}^{t-1} x_{j,h}^{t'} ) . p_{j,h} . T^{a_{j,h}} \right)$. The following cases exist when analyzing Relation (42).

*Case 1.* If all of the remaining work of job *j* is processed within period *t*, then:

*Case 1.1.* If there is no operation $O_{k,l}$ ($k \neq j$) that $z_{k,l,j,h}^{t} = 1$ for each operation $O_{j,h}$ ($h \leq h_j$), then $Tr_{j,t}$ obtained by Relation (42) is equal to $Tr_{j,t}$ obtained by Relation (21).

*Case 1.2.* If at least there is one operation $O_{k,l}$ ($k \neq j$) that $z_{k,l,j,h}^{t} = 1$ for operation $O_{j,h}$ ($h \leq h_j$), then $Tr_{j,t}$ obtained by Relation (42) is smaller than the $Tr_{j,t}$ obtained by Relation (21).

*Case 2.* If part of the remaining work of job *j* is processed in period $t'$ ($t' > t$), then it is evident that $Tr_{j,t}$ obtained by Relation (42) is smaller than $Tr_{j,t}$ obtained by Relation (21).

Combining these two cases shows that $Tr_{j,t}$ obtained by Relation (42) is always smaller than or equal to $Tr_{j,t}$ obtained by Relation (21), and the proof is complete. ☐

Relation (42) balances the remaining work and time for each job. Given this relation, it can be predicted that the model, at each iteration, sets the amount of remaining work for a job based on its remaining time according to its value. It should be noted that in the rolling horizon procedure, if the number of current period $t$ is greater than $tp_j$ ($t > tp_j$), then Constraint (42) is replaced with the following relation.

$$Tr_{j,t} \geq \left( \sum_{h=1}^{h_j} ( \sum_{t'=1}^{t''} D_{j,t'} - \sum_{t'=1}^{t} x_{j,h}^{t'} ) . p_{j,h} . T^{a_{j,h}} \right) - \left( d_{j,t} - (L.(t-1)) \right) . \frac{fr_{j,t}}{k.t} \qquad j = 1,..,J; \; t'' \in tp_j; \; t > t'' \qquad (43)$$

where $fr_{j,t}$ is the remaining lot of job *j* at the end of *the t-th* period, which must be processed in the following production periods (period $t + 1$ to $T$). This parameter can be obtained by Constraint (44). Furthermore, *k* is a constant, which is defined experimentally.

$$fr_{j,t} \geq \sum_{h=1}^{h_j} (D_{j,t} - \sum_{t'=1}^{t} x_{j,h}^{t'}) \qquad j = 1,..,J;\ t \geq tp_j \qquad (44)$$

*5.2. Local search heuristic*

To improve the constructed solutions by the proposed rolling horizon approach, a local search heuristic is developed. The local search is based on maintaining some production decisions through fixed variables while releasing and recalculating others. In other words, all decision variables for a specific job $j$ include $y_{j,h}^{t}$, $z_{j,h,k,l}^{t}$ and $x_{j,h}^{t'}$ ($t = 1,..,T$ ; $k = 1,..,J$; $h = 1,..,h_j$; $l = 1,..,h_k$) are released, and others are fixed. Then, the developed MILP can be solved by commercial MILP solvers. This method is a *job-based local search* and leads to a reduction of overtime and tardiness costs by:

- Moving the entire or part of an item processing from a given period to other periods.
- Altering the position of a given item in the sequence for a specific period.

An essential part of the method is determining which job variables and how many should be released. Maximum improvement in a solution can be obtained in a short time if jobs are properly selected. Even though the computational time for the local search is limited, it significantly influences the total computational time.

Rising the number of released decision variables may result in higher-quality solutions owing to the enlarged search space. Nevertheless, this neighborhood search method is ineffective if many jobs are selected and released. Moreover, this procedure is ineffective because of the choice of a few or the wrong jobs. Several tests have shown that selecting 15 percent of jobs and releasing them is a good alternative. This percentage can be significantly increased for smaller problems. To choose these jobs, first, the tardiness cost of all jobs is calculated separately; next, the jobs are selected based on the following guidelines.

- Release the decision variables for 5% of the jobs with the highest tardiness cost.
- Release the decision variables for 5% of the jobs with the medium tardiness cost.
- Release the decision variables for 5% of the jobs with the lowest tardiness cost.

## 6. Post-processing approach

The material balance represents the most complex part of the multi-level ILSP. In this system, workstations receive materials from lower levels, process them, and then transfer the processed materials to the next level. In the lot-sizing and scheduling problem, it is essential to synchronize the start and finish times of operations within a period with the material balance constraints. The literature offers two primary approaches for establishing material balance in multi-level lot-sizing and scheduling problems. In the first approach, it is assumed that successor items cannot begin production until all predecessor items have been completed. The second approach, known as *lot-streaming*, allows overlapping production between predecessor and successor items. This is achieved by dividing a processing lot of predecessor items into smaller sub-lots, which are gradually transferred to successor operations as they are completed.

In Model-II, lot-streaming is excluded by incorporating the cutting plane (25). However, in some manufacturing environments, lot-streaming is feasible and desirable. To address this, a post-processing approach has been developed to incorporate lot-streaming into the final outputs of the model and algorithms. This post-processing method shifts certain operations to earlier periods, provided this adjustment does not interrupt successor operations due to insufficient material flow from predecessors. This allows successor items to begin before the completion of all predecessor items, enabling overlapping production.

The post-processing approach solves a modified model, where the sequencing variable $z^{t}_{j,h,k,l}$ is fixed to its final value ($\bar{z}^{t}_{j,h,k,l}$) as determined by Model-II or the proposed matheuristic algorithms. The new variables and constraints introduced in the post-processing model are outlined below.

**New variable:**

$\theta^{t,t'}_{j,h}$ Production quantity of operation $O_{j,h}$ is transferred from period $t'$ to $t$ ($t' > t$).

**Post-processing model:**

Constraints (21), (23), (24), (29), (30), (31), (32), (35), (36), (37), (38), (39).

Constraint (34) after fixing $z^{t}_{j,h,k,l}$ to $\bar{z}^{t}_{j,h,k,l}$.

$$x^{t}_{j,h} = \bar{x}^{t}_{j,h} + \sum_{t'=t+1}^{T} \theta^{t,t'}_{j,h} - \sum_{t'=1}^{t-1} \theta^{t',t}_{j,h} \qquad j = 1,..,J; h = 1,..,h_j;\ t = 1,..,T \qquad (45)$$

$$\sum_{t'=1}^{t-1} \theta^{t',t}_{j,h} \leq \bar{x}^{t}_{j,h} \qquad j = 1,..,J; h = 1,..,h_j;\ t = 1,..,T \qquad (46)$$

$$I^{t-1}_{j,h_j} + x^{t}_{j,h_j} - \sum_{t'=1}^{t} \rho^{t',t}_{j} D_{j,t'} = I^{t}_{j,h_j} \qquad j = 1,...,J;\ t = 1,...,T \qquad (47)$$

$$I^{t-1}_{j,h} + x^{t}_{j,h} - x^{t}_{j,h+1} = I^{t}_{j,h} \qquad j = 1,..,J; h = 1,..,h_j - 1;\ t = 1,..,T \qquad (48)$$

$$\rho^{t,t'}_{j} = 0 \qquad j = 1,..,J;\ t = 1,..,T; t' = t,...,T; \text{ If } \rho^{t,t''}_{j} = 1 \text{ where } t'' < t' \qquad (49)$$

$$s^{t}_{j,h+1} \geq s^{t}_{j,h} - \left(I^{t}_{j,h}.\, p_{j,h+1}.\, t^{a_{j,h+1}}\right) \qquad j = 1,..,J; h = 1,..,h_j;\ t = 1,..,T; \text{ If } p_{j,h} \leq p_{j,h+1} \qquad (50)$$

$$s^{t}_{j,h+1} \geq s^{t}_{j,h} + \left(p_{j,h} - p_{j,h+1}\right).\left(x^{t}_{j,h+1} - I^{t-1}_{j,h}\right) + \left(1 - I^{t-1}_{j,h}\right).\, p_{j,h+1}.\, t^{a_{j,h+1}} \qquad j = 1,..,J; h = 1,..,h_j; t = 1,..,T; \text{ If } p_{j,h} > p_{j,h+1} \qquad (51)$$

Constraint (45) determines the updated values of the lot-size variables. In this constraint the post-processing lot-size for operation $O_{j,h}$ in period $t$ is calculated as follows: it equals the pre-processing lot-size determined by Model-II or the proposed matheuristic ($\bar{x}^{t}_{j,h}$), plus the production quantity of $O_{j,h}$ transferred from later periods to period $t$, minus the production quantity of $O_{j,h}$ shifted from period $t$ to earlier periods. Constraint (46) ensures that the production quantity of $O_{j,h}$ transferred from period $t$ to earlier periods cannot exceed the post-processing lot-size value of $O_{j,h}$ in period $t$ ($\bar{x}_{j,h}^{t}$). Constraints (47) and (48) ensure inventory balance for the operations of each job across macro periods. Constraint (49) ensures the post-processing results are only a left-shifted version of the pre-processing results. Constraints (50) and (51) establish precedence relationships between two consecutive operations within the same job while accounting for lot-streaming. Constraint (50) ensures that the successor operation cannot start until the minimum delay required after the predecessor operation has begun, particularly when the predecessor is processed faster than the successor. Conversely, Constraint (51) addresses the same scenario where the predecessor operation progresses slower than the successor.

**Theorem 1.** For the uninterrupted production of $O_{j,h+1}$ in period $t$, the minimum required delay to start $O_{j,h+1}$ after the start time of $O_{j,h}$ is given Constraints (50) and (51).

***Proof.*** The validity of this theorem is detailed in Section 4 of Rohaninejad and Hanzálek [30].

## 7. Computational results

This paper presents the following versions of the rolling horizon heuristic founded on the two mentioned lower bound and local search.

- **RH1**: The rolling horizon heuristic founded on Model-II and lower bound 1.
- **RH2**: The rolling horizon heuristic founded on Model-II and lower bound 2.
- **RH1-LO**: The rolling horizon heuristic founded on Model-II, lower bound 1, and local search.

The proposed MILP models, along with the three heuristic methods, were coded in GAMS 24.8.2, and the CPLEX solver was used to handle the problem-solving tasks. All computations were executed on a personal computer with a 2.66 GHz processor and 6 GB of RAM. The Relative Percentage Deviation (RPD) [45, 46] was employed as the performance metric, which evaluates the deviation of each solution compared to the best solution obtained from Model-I, Model-II, RH1, RH2, and RH1-LO.

$$RPD_A = \frac{(objectve\ function)_A - Min(objectve\ function)}{Min(objectve\ function)} \qquad (52)$$

### *7.1. Generating random instances*

Several random instances are generated to evaluate the performance of the proposed matheuristics. The range for the parameters was adopted from our observation of the case study in the core analysis laboratory. Table 3 indicates how these instances are generated. Operations are randomly allocated to machines. Finally, 15 random instances according to the above conditions are generated, and each instance is labeled with $(\alpha{:}\,\beta{:}\,\gamma{:}\,\delta)$, which respectively represent the number of jobs, operations, machines, and periods within the planning horizon. To establish the best trade-off between the speed and quality of the proposed methods, the maximum computational time of 3600 seconds is determined on the execution of RH1, RH2, and the MILP models. This limit is increased to 7200 seconds for larger instances and reduced to 1800 seconds for smaller instances. In addition, the time limit of 1800 seconds is set for the local search heuristic.

**Table 3.** Values for generating random instances

| Parameter | Notation | Value |
|---|---|---|
| Available Capacity | $C_{i,t}$ | (Beta distribution with $\alpha = 2$ and $\beta = 3$) . $L$ |
| Upper limit of overtime | $O_{m,t}$ | $L - C_{i,t}$ |
| Tardiness costs | $TC_j$ | Uniform distribution between [50 to 200] |
| Overtime costs | $OC_m$ | Uniform distribution between [30 to 60] |
| Processing time rate | $p_{j,h}$ | Beta distribution with $\alpha = 4$ and $\beta = 2$ between [0.5 to 2] |
| Learning index | $a_{j,h}$ | Uniform distribution between [-0. 5 to -0.05] |
| Due date | $d_{j,t}$ | ( Beta distribution with $\alpha = 2$ and $\beta = 4$) . t. $L$ |
| Demands quantity | $D_{j,t}$ | Normal distribution with $\mu = 30$ and $\delta = 8$ |

### *7.2. Performance of solution methods*

To compare Model-I and Model-II, five small- and medium-sized instances are solved by a CPLEX solver with a computational time limit of 1800 seconds, and the obtained results are shown in Table 4. Based on these results, Model-II drastically reduces computational time with an acceptable average RPD. Furthermore, considering Fig. 3, the superiority of Model-II to Model-I, regarding the obtained objective values, by increasing instance size is quite evident. Removing the optimal space of the original problem by applying the cutting plane (Relation 25) caused the obtained solutions of Model-II, for instances TP 2:8:2:3 and TP 3:10:3:3, to be worse than the solutions obtained by Model-I.

**Table 4.** Comparing the performance of Model-I and Model-II

| Instance | Model-I | | | Model-II | | |
|---|---|---|---|---|---|---|
| | Time (s) | Obj. | RPD (%) | Time (s) | Obj. | RPD (%) |
| TP 4:8:2:3 | 1465 | 4928.2 | 0.0% | 1.1 | 5232.5 | 6.2% |
| TP 3:12:4:3 | >3600 | 20187.1 | 0.0% | 1.9 | 20523.6 | 1.7% |
| TP 4:20:5:3 | >-3600 | 27976.6 | 12.5% | 5.8 | 24877.9 | 0.0% |
| TP 5:25:5:3 | >3600 | 37513.1 | 26.3% | 81.3 | 29690.5 | 0.0% |
| TP 6:30:5:3 | >3600 | 55244.4 | 57.57% | >3600 | 35059.5 | 0.0% |
| Average | >3173 | - | 19.3% | >738 | - | 1.6% |

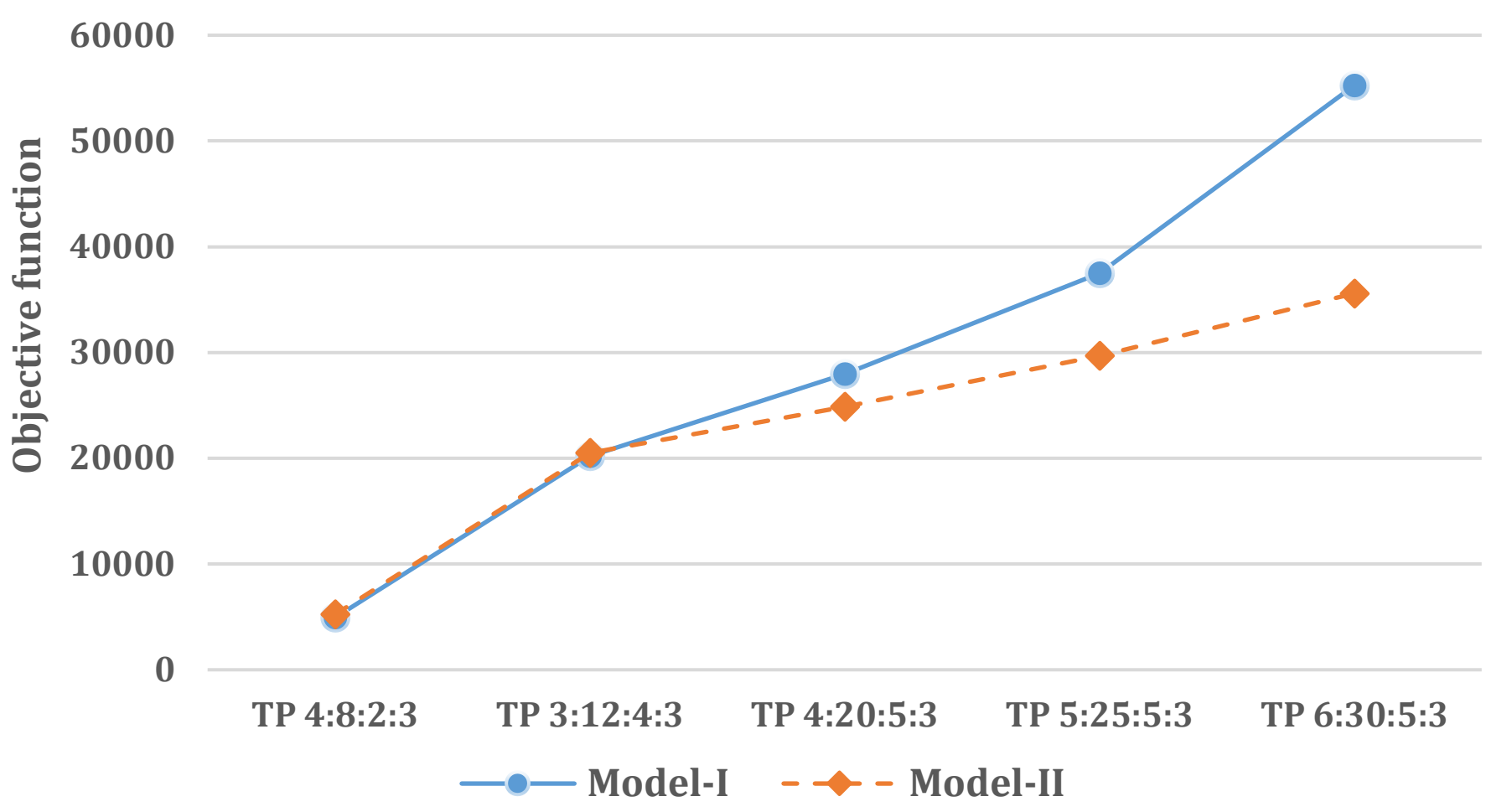

**Fig. 3.** Comparing the objective function of Model-I and Model-II

Table 5 presents the performance of the proposed lower bounds, which were based on Model-II. Regarding this table, GAP 1 and GAP 2 define the percentage of the difference between the optimal solution and lower bounds 1 and 2, respectively. Based on the results, the performance of both lower bounds is satisfying, especially since the objectives of lower bound 1 are considerably close to the optimal values.

**Table 5.** Performance of the proposed lower bounds in Model-II

| Instance | Model-II | | | | Lower bound 1 | | Lower bound 2 | |
|---|---|---|---|---|---|---|---|---|
| | Time (s) | Obj. | GAP 1 (%) | GAP 2 (%) | Time (s) | Obj. | Time (s) | Obj. |
| TP 2:8:2:3 | 1.1 | 15102.3 | 2.6 | 18.6 | 1.6 | 14717.9 | 1.1 | 12736.5 |
| TP 3:10:3:3 | 1.8 | 23789.8 | 3.9 | 25.8 | 1.2 | 22892.2 | 1.3 | 18917.6 |
| TP 4:15:4:3 | 3.2 | 33029.2 | 2.3 | 26 | 1.6 | 32296.1 | 1.2 | 26221.8 |
| TP 5:20:5:3 | 6.4 | 40240.7 | 3.8 | 19.3 | 3.5 | 38752.4 | 1.4 | 33725.7 |
| TP 5:25:5:3 | 11.7 | 47132.9 | 6.3 | 15.8 | 5.4 | 44319.3 | 1.2 | 40712.9 |
| Average | 4.8 | - | 3.7 | 21.1 | 2.6 | - | 1.4 | - |

To assess and contrast the behavior of the extended solution methods, 10 large- and very large-sized instances are solved, and the obtained results are juxtaposed with the results of Model-II in Table 6. Based on the size of the problems, computational time limits of 3600 and 7200 seconds are considered. This table shows that Model-II cannot even reach a feasible solution for large-sized instances in 7200 seconds. According to the RPD values, the quality of the proposed heuristics is clearly satisfying. In addition, the RH1-LO heuristic method surpasses other methods remarkably when the size of problems increases. Although this fact was expected somehow, the improvement in the results of RH1 by applying the proposed local search is remarkable. The objective values obtained by RH1 are less than those obtained by RH2 and Model-II. Therefore, combining lower bound 1 with a rolling horizon procedure has a more positive impact than using lower bound 2. However, in terms of computational time, RH2 is the superior method. Moreover, as Fig. 4 shows, the difference between the quality of RH1 and RH2 decreases in very large-sized instances so that the objective function value for the RH2 is less than RH1 in (TP 50:500:22:5).

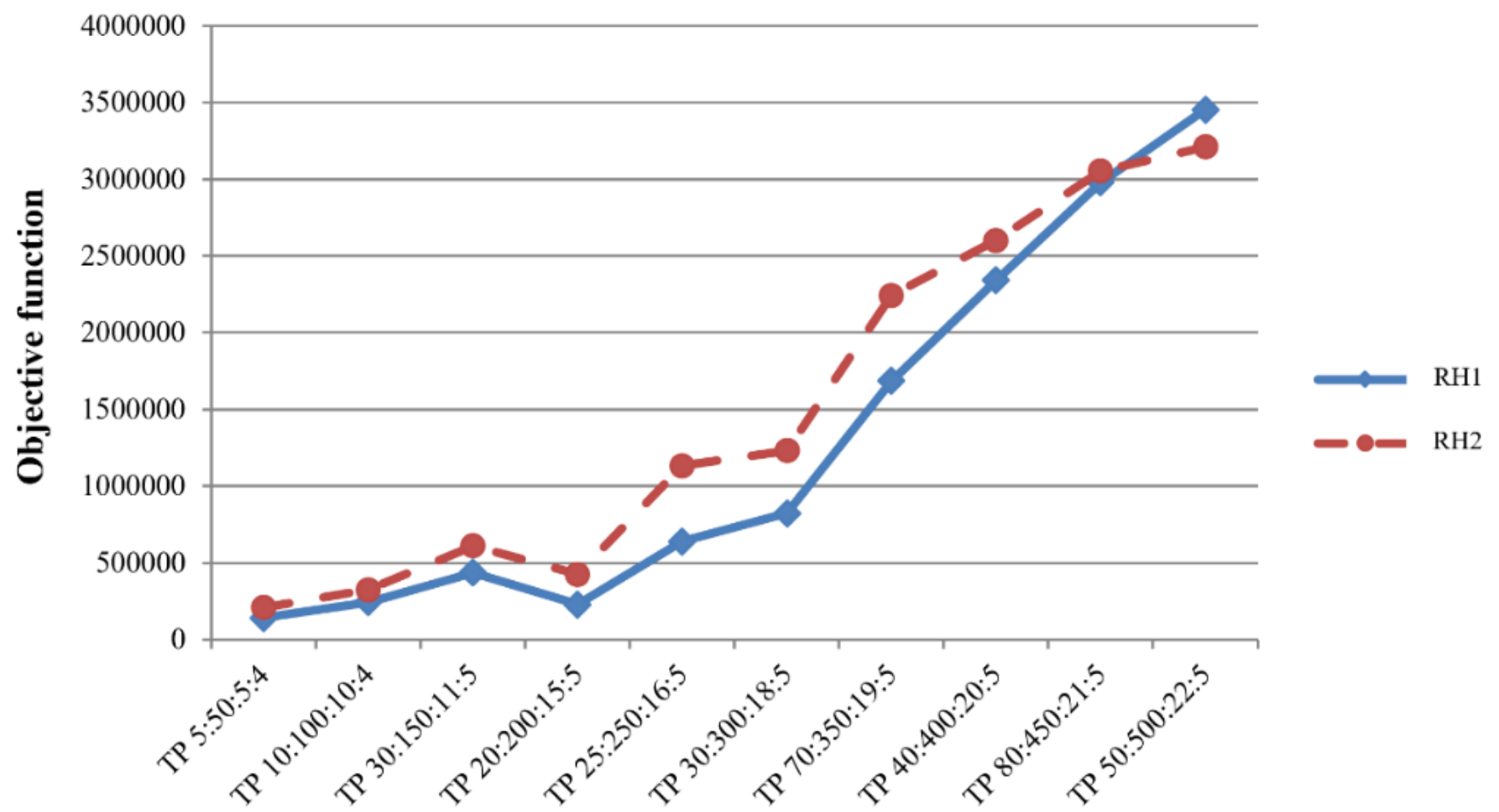

**Fig. 4.** Comparing the objective function of RH1 and RH2

Tables 7 to 10 exhibit the impact of the learning index on the performance of the proposed lower bounds. GAP 1 and GAP 2 are defined as the percentage of the difference between the optimal solutions and the solutions gained by lower bounds 1 and 2, respectively. According to the data provided in these tables, it can be inferred that the performance of lower bound 1 is independent of the learning effect. In contrast, the performance of lower bound 2 is improved by reducing the learning effect. Hence, lower bound 2 for problems without the learning effect can be strongly advisable.

It should be noted that increasing the learning effect has a negative impact on lower bound 2. This negative effect is due to the difference between $tp_j$ and $T$. This reason can predict the decline in the performance of lower bound 2 with increasing the number of periods ($T$) as well. However, the following relation can be a substitute for Relation (42) in problems with a high learning effect or many periods. Although it cannot be proven that this relation is a lower bound for Model-II, it can experimentally create more favorable results.

**Table 6.** Comparing the performance of the developed heuristics and Model-II

| Instance | Model-II | | | RH1 | | | RH2 | | | RH1-LO | | |
|---|---|---|---|---|---|---|---|---|---|---|---|---|
| | Time (s) | Obj | RPD (%) | Time (s) | Obj | RPD (%) | Time (s) | Obj | RPD (%) | Time (s) | Obj | RPD (%) |
| TP 5:50:5:4 | >3600 | 120306 | 0.0 | 12.1 | 143602 | 19.4 | 3.6 | 210620 | 75.1 | 25.7 | 137550 | 14.3 |
| TP 10:100:10:4 | >3600 | 213589 | 0.0 | 39.3 | 245870 | 15.1 | 14.2 | 325441 | 52.4 | 82.3 | 231747 | 8.5 |
| TP 30:150:11:5 | >3600 | 649991 | 56.7 | 3335 | 436163 | 5.2 | 258 | 612859 | 47.8 | 3840 | 414791 | 0.0 |
| TP 20:200:15:5 | >3600 | 767136 | 241.5 | 368 | 231585 | 3.1 | 132 | 425510 | 89.5 | 522 | 224568 | 0.0 |
| TP 25:250:16:5 | >7200 | 1604625 | 159.1 | 7128 | 640528 | 3.4 | 1283 | 1135600 | 83.4 | >8928 | 619230 | 0.0 |
| TP 30:300:18:5 | >7200 | 3248552 | 311.8 | 5159 | 823741 | 4.4 | 882 | 1232068 | 56.2 | >6959 | 788922 | 0.0 |
| TP 70:350:19:5 | >7200 | * | * | >7200 | 1687323 | 4.4 | 3955 | 2242752 | 38.8 | >9000 | 1615749 | 0.0 |
| TP 40:400:20:5 | >7200 | * | * | 2425 | 2342066 | 3.2 | 326 | 2602337 | 14.7 | 3407 | 2268845 | 0.0 |
| TP 80:450:21:5 | >7200 | * | * | >7200 | 2977462 | 5.1 | 4762 | 3051365 | 7.7 | >9000 | 2833400 | 0.0 |
| TP 50:500:22:5 | >7200 | * | * | 3856 | 3452410 | 7.6 | 418 | 3209490 | 0.0 | 5192 | 3248637 | 1.2 |
| Average | >7200 | — | 128.2 | >3672 | — | 7.1 | 1203 | — | 46.5 | >4695 | — | 2.4 |

$$Tr_{j,t} \geq \left( \sum_{h=1}^{h_j} \left( \sum_{t'=1}^{t} D_{j,t'} - \sum_{t'=1}^{t-1} x_{j,h}^{t'} \right) . p_{j,h} . t^{a_{j,h}} \right) - \left( d_{j,t} - \left( L.(t-1) \right) \right) \qquad j = 1,..,J; \; t \in tp_j \tag{53}$$

**Table 7.** Performance of the proposed lower bounds for Model-II without the learning effect

| Instance | Model-II | | | | Lower bound 1 | | Lower bound 2 | |
|---|---|---|---|---|---|---|---|---|
| | Time (s) | Obj. | GAP 1 (%) | GAP 2 (%) | Time (s) | Obj. | Time (s) | Obj. |
| TP 2:8:2:3 | 1.1 | 18684.0 | 3.2 | 13.4 | 1.1 | 18100.3 | 1.1 | 16472.1 |
| TP 3:10:3:3 | 1.2 | 30728.7 | 4.6 | 19.0 | 1.1 | 29364.1 | 1.1 | 25822.6 |
| TP 4:15:4:3 | 2.3 | 40687.4 | 6.2 | 14.9 | 1.7 | 38313.1 | 1.2 | 35417.2 |
| TP 5:20:5:3 | 12.4 | 53822.4 | 7.7 | 13.4 | 3.2 | 49964.2 | 1.3 | 47452.7 |
| TP 5:25:5:3 | 49.2 | 74186.5 | 7.6 | 15.8 | 11.4 | 68929.1 | 1.6 | 64070.4 |
| Average | 13.2 | - | 5.9 | 15.3 | 3.7 | - | 1.2 | - |

**Table 8.** Performance of the proposed lower bounds for Model-II with the learning index - 0.2

| Instance | Model-II | | | | Lower bound 1 | | Lower bound 2 | |
|---|---|---|---|---|---|---|---|---|
| | Time (s) | Obj. | GAP 1 (%) | GAP 2 (%) | Time (s) | Obj. | Time (s) | Obj. |
| TP 2:8:2:3 | 1.1 | 15880.5 | 0.3 | 31.5 | 1.1 | 15832.4 | 1.0 | 12075.9 |
| TP 3:10:3:3 | 1.1 | 26073.4 | 4.1 | 23.7 | 1.1 | 25046.3 | 1.0 | 21078.9 |
| TP 4:15:4:3 | 2.9 | 35377.3 | 3.0 | 39.1 | 1.8 | 34362.3 | 1.0 | 25426.2 |
| TP 5:20:5:3 | 6.2 | 44951.2 | 3.9 | 27.8 | 2.2 | 43277.4 | 1.4 | 35168.1 |
| TP 5:25:5:3 | 17.5 | 57484.8 | 5.8 | 16.4 | 4.4 | 54356.9 | 3.0 | 49390.9 |
| Average | 5.7 | - | 3.4 | 27.7 | 2.1 | - | 1.4 | - |

**Table 9.** Performance of the proposed lower bounds for Model-II with the learning index - 0.4

| Instance | Model-II | | | | Lower bound 1 | | Lower bound 2 | |
|---|---|---|---|---|---|---|---|---|
| | Time (s) | Obj. | GAP 1 (%) | GAP 2 (%) | Time (s) | Obj. | Time (s) | Obj. |
| TP 2:8:2:3 | 1.1 | 13392.3 | 0.0 | 58.8 | 1.1 | 13392.3 | 1.1 | 8435.6 |
| TP 3:10:3:3 | 1.2 | 22573.3 | 2.1 | 37.7 | 1.1 | 22103.7 | 1.1 | 16392.0 |
| TP 4:15:4:3 | 2.8 | 31066.7 | 5.8 | 65.4 | 2.8 | 29351.0 | 1.1 | 18782.2 |
| TP 5:20:5:3 | 5.3 | 38442.7 | 3.8 | 44.4 | 1.8 | 37022.9 | 1.2 | 26631.3 |
| TP 5:25:5:3 | 7.9 | 48072.5 | 4.1 | 34.0 | 3.9 | 46185.7 | 1.5 | 35868.0 |
| Average | 3.6 | - | 3.2 | 48.1 | 2.1 | - | 1.2 | - |

**Table 10.** Performance of the proposed lower bounds for Model-II with the learning index - 0.6

| Instance | Model-II | | | | Lower bound 1 | | Lower bound 2 | |
|---|---|---|---|---|---|---|---|---|
| | Time (s) | Obj. | GAP 1 (%) | GAP 2 (%) | Time (s) | Obj. | Time (s) | Obj. |
| TP 2:8:2:3 | 1.1 | 11328.1 | 0.0 | 148 | 1.1 | 11327.7 | 1.0 | 4564.5 |
| TP 3:10:3:3 | 1.0 | 19682.7 | 6.2 | 94 | 1.0 | 18540.6 | 0.9 | 10123.4 |
| TP 4:15:4:3 | 1.9 | 26347.9 | 7.9 | 126 | 1.3 | 24418.7 | 1.0 | 11673.8 |
| TP 5:20:5:3 | 4.1 | 33754.8 | 2.4 | 83 | 1.5 | 32968.5 | 1.0 | 18486.2 |
| TP 5:25:5:3 | 6.0 | 40170.2 | 2.6 | 69 | 1.6 | 39146.0 | 1.1 | 23742.4 |
| Average | 2.8 | - | 3.8 | 104 | 1.6 | - | 1.1 | - |

Table 11 compares the results of Model-II before and after incorporating lot-streaming. As outlined in Section 6, the post-processing approach facilitates overlapping production between predecessor and successor operations by progressively transferring smaller sub-lots from the predecessor to the successor. The results demonstrate that incorporating lot-streaming improves the objective function by an average of 58%, which is a significant enhancement. Furthermore, the post-processing approach adds minimal computational overhead to the solution process. This combination of substantial improvement with low additional complexity makes the post-processing approach a practical and effective solution for manufacturing systems where lot-streaming is possible and desirable.

**Table 11.** Model-II comparison: Pre- and Post-Processing with Lot-Streaming

| Instance | Pre-Processing (Without Lot-Streaming) | | | Post-Processing (With Lot-Streaming) | | | |
|---|---|---|---|---|---|---|---|
| | Time (s) | Obj. | RPD (%) | Time (s) | Total Times (s) | Obj. | RPD (%) |
| TP 4:8:2:3 | 1.1 | 5232.5 | 5.2% | 0.1 | 1.2 | 4972.4 | 0.0% |
| TP 3:12:4:3 | 1.9 | 20523.6 | 78.8% | 0.17 | 2.07 | 11480.8 | 0.0% |
| TP 4:20:5:3 | 5.8 | 24877.9 | 72.8% | 0.15 | 5.95 | 14399.4 | 0.0% |
| TP 5:25:5:3 | 81.3 | 29690.5 | 59.2% | 0.16 | 81.46 | 18652.3 | 0.0% |
| TP 6:30:5:3 | >3600 | 35059.5 | 75.5% | 0.31 | >3600 | 19972.8 | 0.0% |
| Average | >738 | - | 58.3% | 0.16 | >738 | - | 0.0% |

To validate the proposed approach, a real case is considered in the context of the Special Core Analysis Laboratory (SCAL), which focuses on advanced rock and reservoir property evaluations. This case involves 6 jobs, 9 operations, and 6 periods, with demand measured in sample units. Related data are available at https://github.com/Rohanmoh/ILSP-Learning-Effect. Priorities between operations arise due to dependencies in experimental processes, which can vary across jobs based on their unique requirements. Lot-streaming was incorporated after a post-processing approach, resulting in an objective function of 12,414. In contrast, the pre-processing stage without lot-streaming produced an objective function of 45,513, obtained using the RH1-LO algorithm in 340 seconds. The resulting Gantt chart, presented in Fig. 5, depicts the best-obtained schedule, with operations for the same job distinguished by matching colors, demonstrating the effectiveness of the applied methodology.

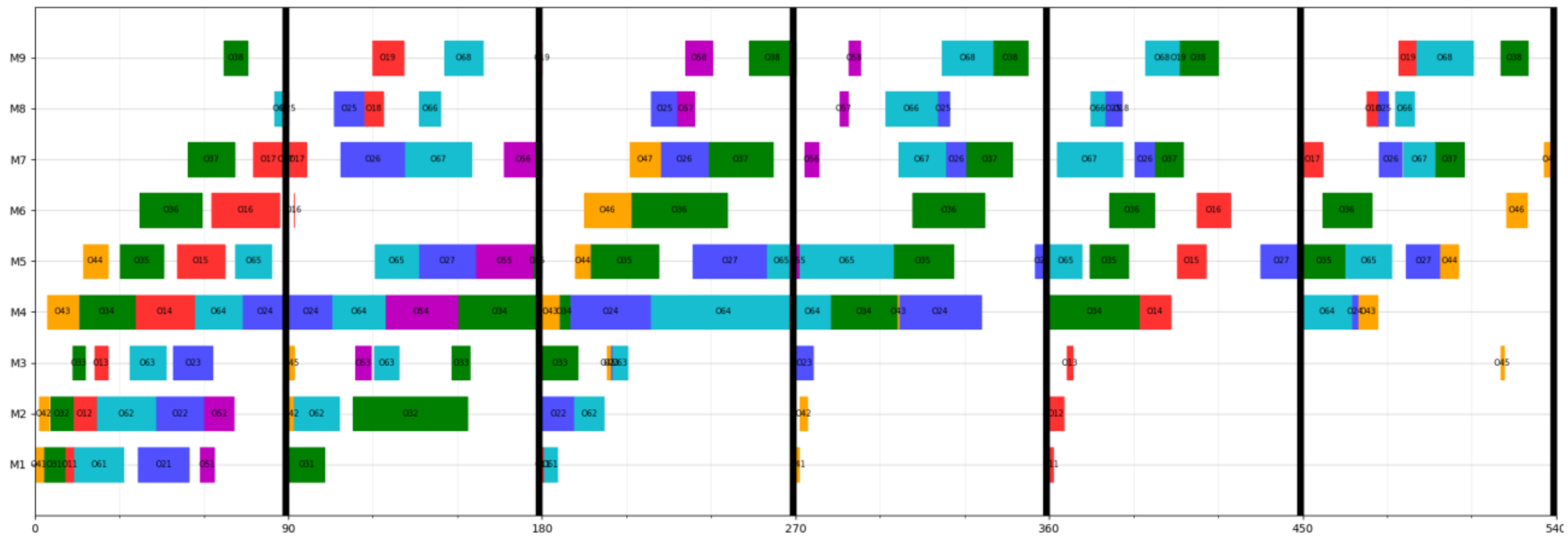

**Fig. 5.** The Gantt chart of the case study with lot-steaming

*7.3. Sensitivity analysis*

In this section, we conduct a sensitivity analysis to examine how variations in key parameters influence the model's performance, enabling better-informed decision-making. The analysis focuses on five critical parameters: processing time, due date, learning rate, tardiness cost, and overtime cost. Each parameter is changed systematically from -30% to +30%, covering seven levels [-30%, -15%, -5%, 0%, +5%, +15%, +30%]. Two sample cases, with dimensions TP 5:30:5:6, are analyzed, and the average percentage of changes in the objective function are depicted in Figs. 6 to 9.

The results reveal that the model's sensitivity is highest for the due date, followed by processing time and the learning rate, with cost-based parameters (tardiness and overtime costs) having a comparatively lower impact. Specifically, increasing the due dates and processing times substantially impacted on the objective function, with maximum variations reaching 2400% and 977%, respectively, highlighting their significant influence. Although the learning effect's influence is periodic, it remains significant, demonstrating its importance even when its impact emerges gradually over a longer planning horizon. Overall, time-based parameters, including due date, processing time, and learning rate, significantly outweigh cost-based parameters in their effect on the model's performance. However, the impact of tardiness and overtime costs, with maximum variations ranging from approximately -18% to +20% in the objective function, remains noteworthy and cannot be overlooked.

These findings provide some managerial implications. Given the pronounced effect of due dates, managers should prioritize strategies to negotiate flexible deadlines with clients, offer discounts or incentives for extended due dates, or explore outsourcing options to meet tight schedules. For processing time, efforts should prevent its increases through proactive maintenance, equipment calibration, staff training, and investing in faster and more reliable machinery. Despite its periodic nature, the learning effect also demands attention; fostering knowledge sharing, maintaining detailed process documentation, and implementing mentorship programs can amplify its benefits. Organizations can effectively enhance their capacity to manage complex lot-sizing and scheduling problems by addressing these areas.

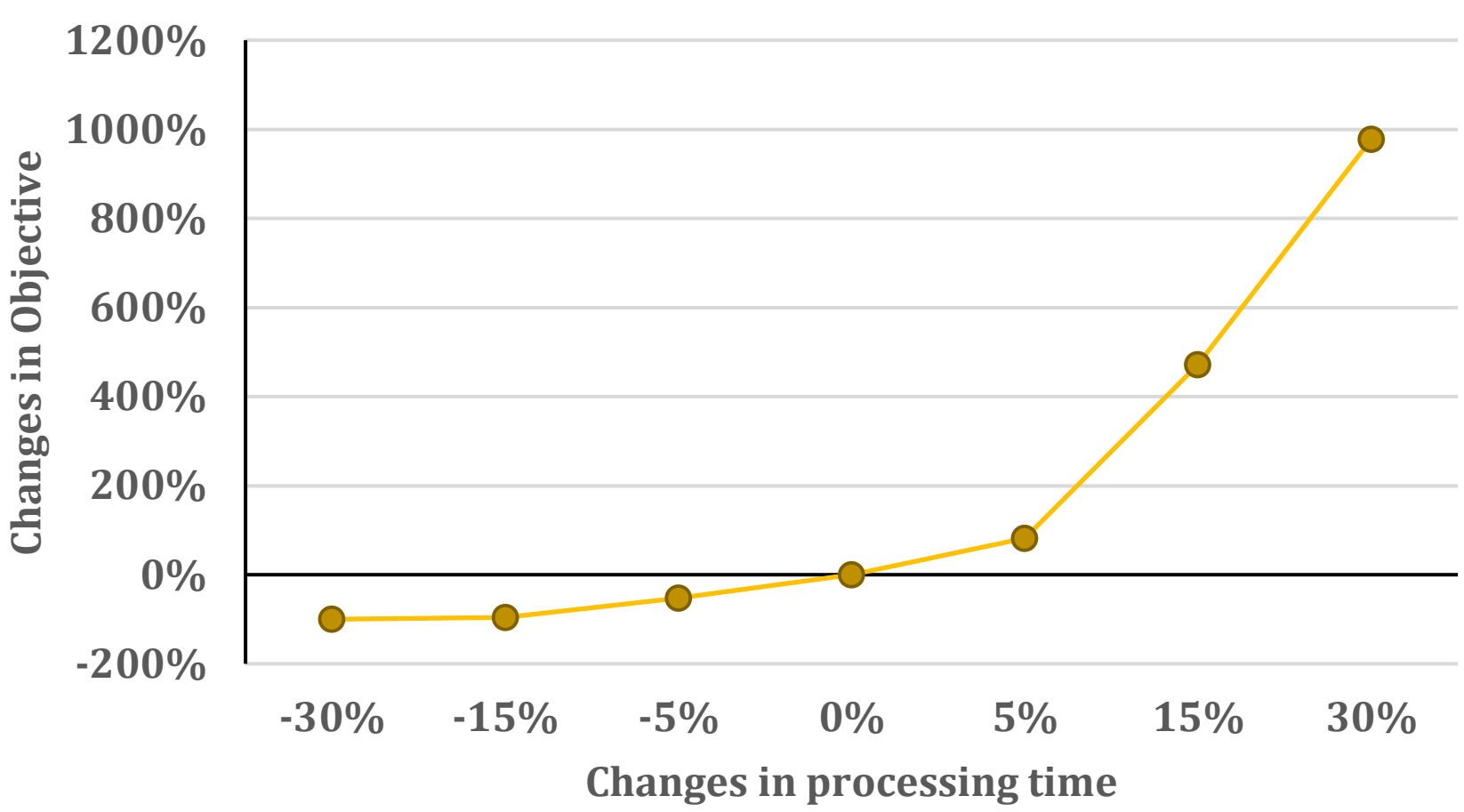

**Fig. 6.** Sensitivity analysis based on the processing time changes

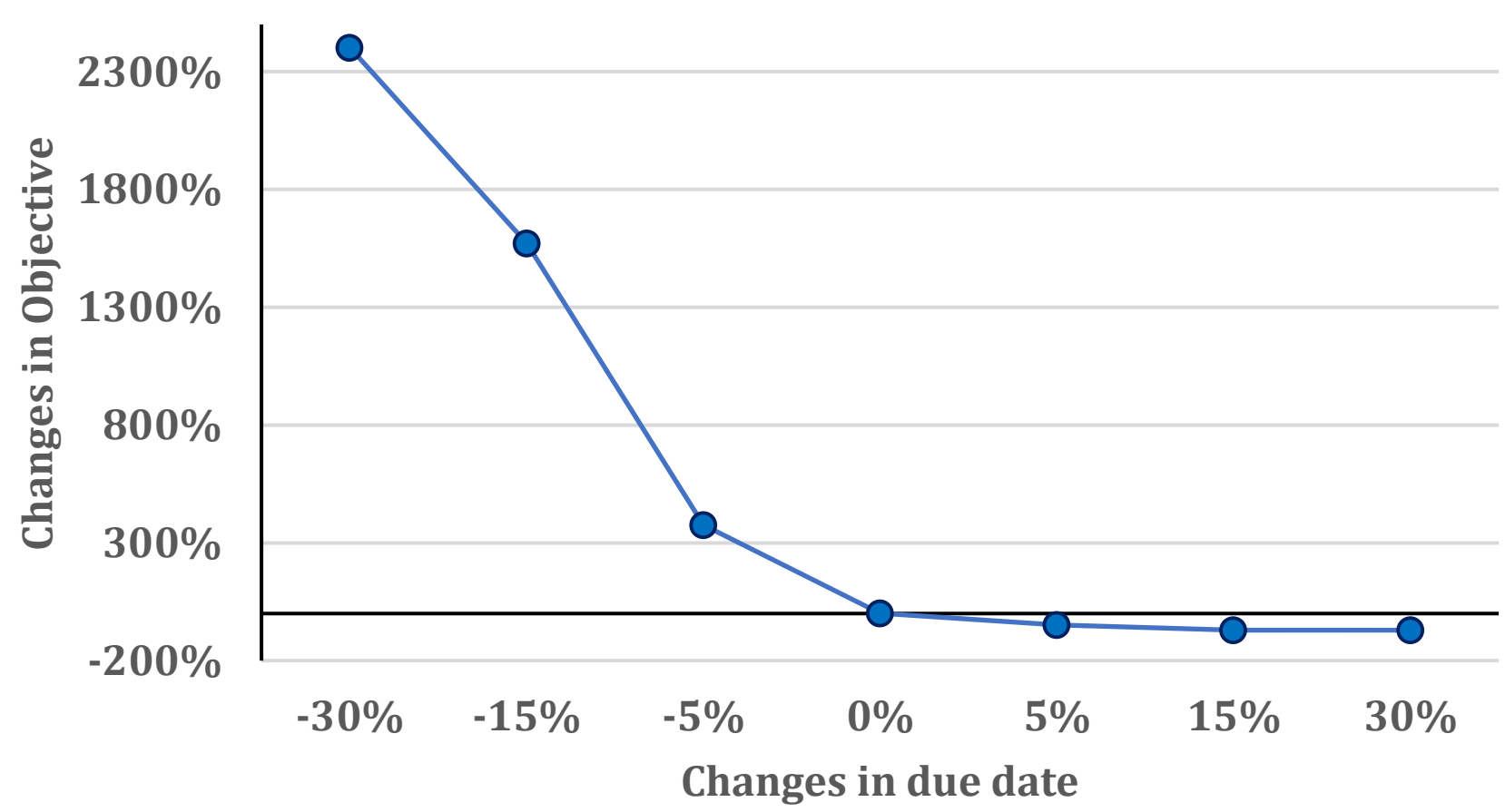

**Fig. 7.** Sensitivity analysis based on the due date changes

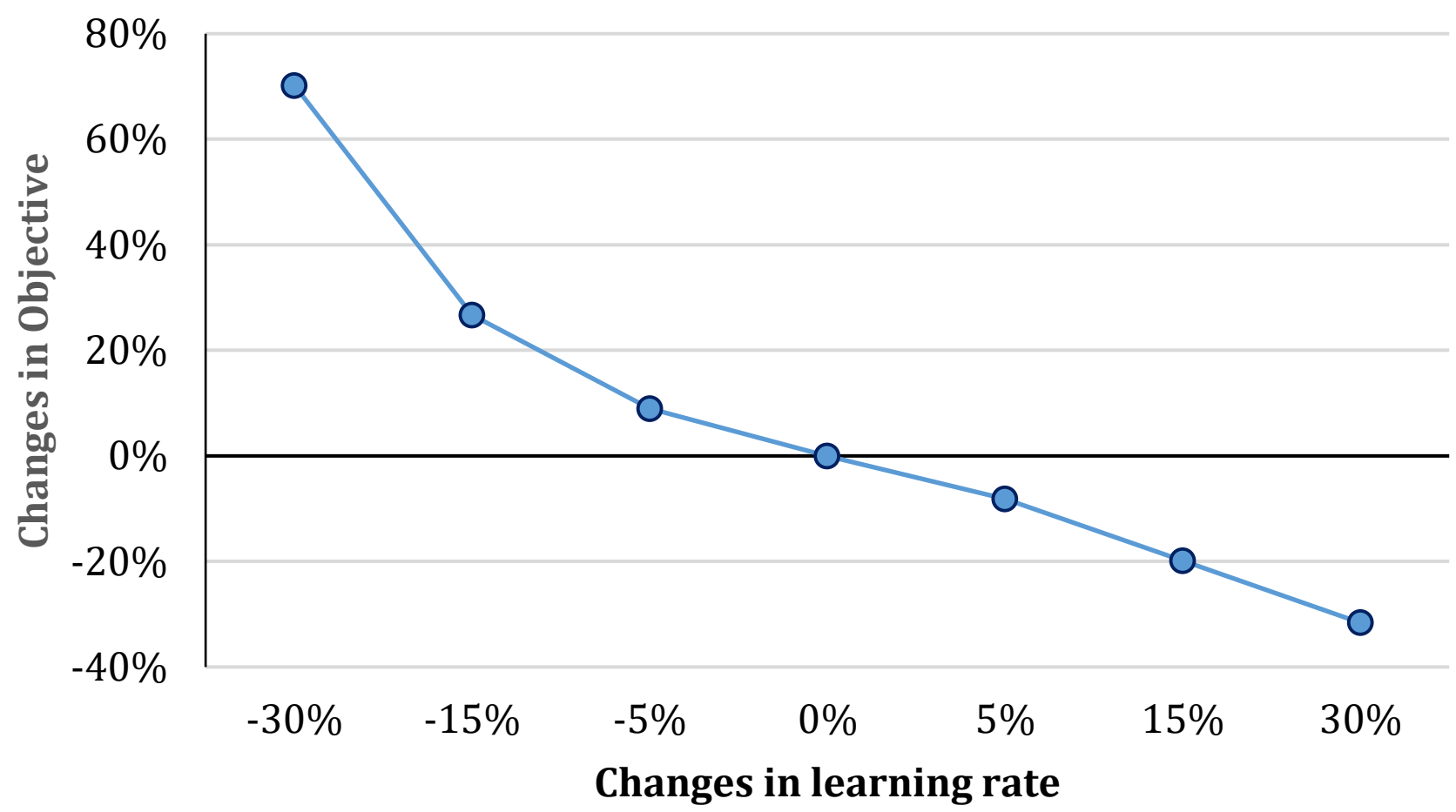

**Fig. 8.** Sensitivity analysis based on the learning effect changes

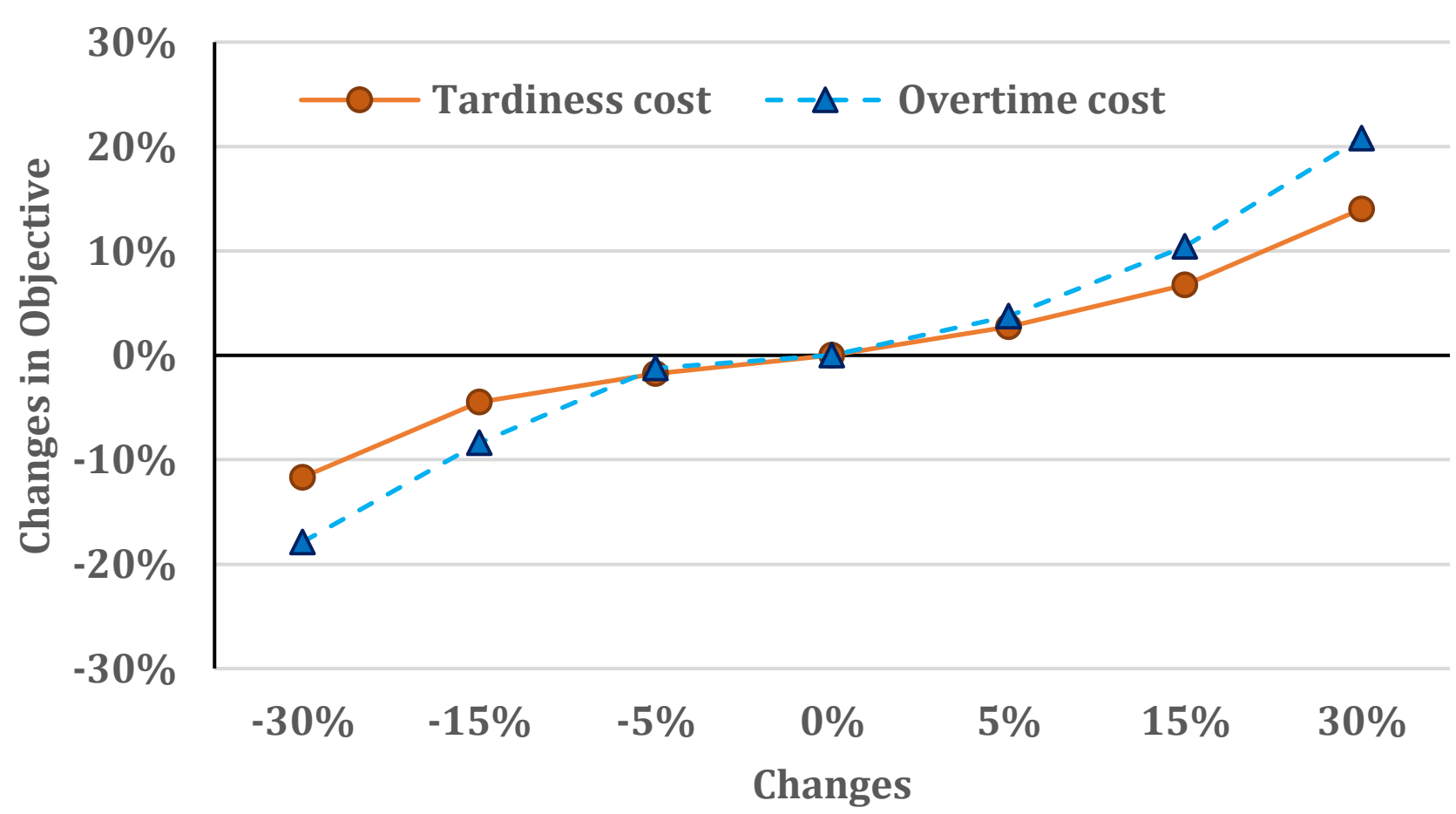

**Fig. 9.** Sensitivity analysis based on the tardiness and overtime cost changes

## 8. Conclusion and future research directions

This study discussed ILSP-JSP with the learning effect. First, a novel big-bucket MILP model was presented and simplified by a cutting plane. Moreover, three matheuristic methods, called RH1, RH2, and RH1-LO, were developed using the rolling horizon approach to solve the problem. In addition, a post-processing approach was employed to incorporate lot-streaming in possible manufacturing settings. Based on the computational results, all three matheuristic methods exhibited superior performance compared to the MILP models, and the RH1-LO surpasses the other heuristics. In addition, the local search method, applied in RH1-LO, improved the results of the RH1 heuristic method by about 5 percent on average. Although the quality of RH1 generally excelled RH2, but RH2 had a lower complexity. Therefore, applying RH2 to very large problems was promising and recommendable. Besides, the results showed that the performance of RH1 was independent of the learning effect, while the performance of RH2 could be improved in cases with lower values of the learning effect. Furthermore, considering the low computational burden of the post-processing approach, it can drastically improve the objective where the lot-streaming is applicable.

According to the sensitivity analysis results, the due date was the most influential parameter on the objective, followed by processing time, learning rate, tardiness cost, and overtime cost. As a result, managers should establish strategies to negotiate flexible due dates, enhance maintenance and renovation practices, and invest in knowledge transfer, mentorship, and training programs for workers to strengthen the performance of their systems.

Extending the solution methods by combining them with meta-heuristic algorithms is suggested so that the binary variables of the MILP models are obtained by a meta-heuristic algorithm and become fixed. Then, a Linear Programming (LP) model solver solves the remaining LP model. Moreover, considering uncertainty in the problem can lead to more realistic results. Furthermore, according to international concerns for tackling global warming and reducing greenhouse gas emissions, incorporating related objectives to the problem can be a valuable direction for future studies.

## CRediT authorship contribution statement

**M. Rohaninejad**: Writing- original draft, Methodology, Software. **B. Vahedi-Nouri**: Conceptualization, Writing - original draft, Validation. **R. Tavakkoli-Moghaddm**: Supervision, Writing – review and editing. **Z. Hanzálek**: Project administration, Writing – review and editing.

## Data Availability Statement

The data supporting this study's findings are available from the corresponding author upon reasonable request.

## Declaration of Competing Interest

The authors declare that they have no known competing financial interests or personal relationships that could have appeared to influence the work reported in this paper.

## Acknowledgment

This work was co-funded by the European Union under the project ROBOPROX (reg. no. CZ.02.01.01/00/22_008/0004590) and the Grant Agency of the Czech Republic under the Project GACR 22-31670S.